%
%
%
%
%
\font\tenbf=cmbx10

\font\eightrm=cmr8
\font\eightit=cmti8

\def\sectiontitle#1\par{\vskip0pt plus.1\vsize\penalty-250
\vskip0pt plus-.1\vsize\bigskip\vskip\parskip
\message{#1}\leftline{\tenbf#1}\nobreak\vglue 5pt}

\def\frac#1#2{{#1\over#2}}

\def\m@th{\mathsurround=0pt}

\def\fsquare(#1,#2){
\hbox{\vrule$\hskip-0.4pt\vcenter to #1{\normalbaselines\m@th
\hrule\vfil\hbox to #1{\hfill$\scriptstyle #2$\hfill}\vfil\hrule}$\hskip-0.4pt
\vrule}}
\magnification=\magstep1
\vsize=8.6truein 
\parindent=15pt
\nopagenumbers

\global\def\ssectitle#1\par{\bigbreak\medskip
  \leftline{\typc #1}
  \nobreak\bigskip\vskip-\parskip\noindent}
\global\def\sssectitle#1\par{\bigbreak\medskip
  \leftline{\typd #1}
  \nobreak\bigskip\vskip-\parskip\noindent}

\def\bitem#1{\item{{[#1]}\quad}}

\font\typc=cmbx10 scaled \magstep1  
\font\typd=cmbx10 scaled \magstep2  
 %
\input amssym.def

\baselineskip=13pt
\headline{\ifnum\pageno=1\hfil\else%
{\ifodd\pageno\rightheadline \else \leftheadline\fi}\fi}
\def\rightheadline{\hfil\eightit 
Domino tableaux, Sch\"utzenberger involution, and the symmetric group action
\quad\eightrm\folio}
\def\leftheadline{\eightrm\folio\quad 
\eightit 
A. Berenstein and A.N. Kirillov 
\hfil}
\voffset=2\baselineskip 
\centerline{\tenbf DOMINO \hskip 0.1cm TABLEAUX, \hskip 0.1cm SCH\"UTZENBERGER
\hskip 0.1cm INVOLUTION,}
\vskip 0.1cm 
\centerline{\tenbf AND \hskip 0.1cm THE \hskip 0.1cm SYMMETRIC \hskip 0.1cm 
GROUP \hskip 0.1cm ACTION
}
\vglue 16pt
\centerline{\eightrm 
ARKADY BERENSTEIN
}
\baselineskip=12pt
\centerline{\eightit 
Department of Mathematics, Cornell University
}
\baselineskip=10pt
\centerline{\eightit 
Ithaca, NY 14853, U.S.A.
}
\vglue 16pt
\centerline{\eightrm 
ANATOL N. KIRILLOV
}
\baselineskip=12pt
\centerline{\eightit 
CRM, University of Montreal
}
\baselineskip=10pt
\centerline{\eightit 
C.P. 6128, Succursale A, Montreal (Quebec) H3C 3J7, Canada
}
\baselineskip=12pt
\centerline{\eightit
and }
\baselineskip=12pt
\centerline{\eightit 
Steklov Mathematical Institute,
}
\baselineskip=10pt
\centerline{\eightit 
Fontanka 27, St.Petersburg 191011, Russia
}
\vglue 20pt

\centerline{\eightrm ABSTRACT}
{\rightskip=1.5pc
\leftskip=1.5pc
\eightrm\parindent=1pc
We define an action of the symmetric group $S_{[{n\over 2}]}$ on the 
set of domino tableaux, and prove that the number of domino tableaux of
weight $\beta'$ does not depend on the permutation of the weight $\beta'$.
A bijective proof of the well-known result due to J.~Stembridge that the
number of self--evacuating tableaux of a given shape and weight $\beta=
(\beta_1,\ldots ,\beta_{[{n+1\over 2}]},\beta_{[{n\over 2}]},\ldots ,
\beta_1)$, is equal to that of domino tableaux of the same shape and
weight $\beta'=(\beta_1,\ldots ,\beta_{[{n+1\over 2}]})$ is given.
} 
\vglue12pt
\baselineskip=13pt
\overfullrule=1pt

\sssectitle {0. Introduction} 

Domino tableaux play a significant role in the representation theory of the 
symmetric group $S_n$, and the theory of symmetric polynomials. The shortest 
definition of semi-standard domino tableaux is given in [10], p. 139. In 
particular, the generating function of the number of domino tableaux of  shape 
$\lambda=(\lambda_1\ge \lambda_2\ge \cdots \ge \lambda_n\ge 0)$, and weight 
$\beta'=(\beta_1,\ldots,\beta_{[ {n+1\over 2}]})$ is given by specializing 
the Schur polynomial $s_\lambda(x)$ at  $x\mapsto y$ where $y_n=-y_{n-1}$, 
$y_{n-2}=-y_{n-3},\ldots$:
$$\pm s_\lambda(y)=\sum_{\beta'}  K^{(2)}_{\lambda,\beta'} \ 
y_n^{2\beta_1}y_{n-2}^{2\beta_2}\cdots \ \ ,\eqno (0.1)$$
where the sign $\pm$ depends only on  $\lambda$. 

Recently they were studied in connection with  the representation theory  of 
$GL_n$.  In particular,  in  [12] John Stembridge proved the following 
conjecture by Richard Stanley (who stated it for $\beta'=(1,1,\ldots,1)$, 
that is, for {\it standard} tableaux).   

\proclaim Theorem 0.1 ([12]). The number of self-evacuating tableaux of shape 
$\lambda$ and  weight $\beta=(\beta_1,\ldots,\beta_{[{n+1\over 2}]},
\beta_{[{n\over 2}]},\ldots,\beta_1)$, equals 
the number of domino tableaux $K^{(2)}_{\lambda,\beta'}$ of  
shape $\lambda$ and weight $\beta'$.

The {\it  self-evacuating tableaux} are defined  as
fixed points of the {\it Sch\"utzenberger evacuation involution} ${\bf S}$
on the {\it semi-standard tableaux} 
(see [11], [12]). The proof of Theorem 0.1 in [12] uses properties 
of the canonical basis
for the $GL_n$-module $V_\lambda$. A key ingredient of the
proof is the result of [3], Theorem 8.2,  which states that ${\bf S}$ is
naturally identified with
the longest element of the symmetric group $S_n$ acting on the canonical 
basis of $V_\lambda$. 
 
There are combinatorial algorithms (see e.g., [5] or [9]) 
which lead to a bijective proof of Theorem 0.1 for 
the {\it standard} tableaux.  For semi-standard tableaux and even $n$ 
a similar algorithm appeared  in [7], remark on the page 399. 

In the present paper we give a bijective proof of Theorem~0.1 for any $n$ 
(Theorem~1.2). 

The main idea of the bijection constructed in this paper is that 
(semi-standard) domino  tableaux are 
naturally identified with those ordinary {\it (semi-standard) tableaux}  
which are fixed under a certain involution ${\bf D}$ (see the Appendix). 

We discover that the involutions ${\bf D}$ and ${\bf S}$ are conjugate as 
automorphisms of all tableaux, say 
${\bf S}={\bf P}{\bf D}{\bf P}^{-1}$
(see Lemma 1.3). Thus, the conjugating automorphism  ${\bf P}$ bijects the 
domino tableaux onto the  self-evacuating tableaux.  

The main difference of our construction and those given in [5], [9]
and [7] is that instead of a description of combinatorial algorithms
as it was done in the papers mentioned, we give a direct algebraic formula 
for the bijection under consideration. Our approach is based on results 
from [2].

It follows from (0.1)  that the number $K^{(2)}_{\lambda,\beta'}$  does not 
depend on the permutation  of the weight $\beta'$. This is analogous to the 
property of the ordinary Kostka numbers $K_{\lambda,\beta}$. 
Recently,
Carre and Leclerc constructed in [4] an appropriate action of the 
symmetric group $S_{n\over 2}$  (for even $n$) which realizes this property. 
Their construction used
a generalized RSK-type algorithm ([4], Algorithm 7.1, and Theorem 7.8). 

The second result of the present paper is an algebraic 
construction of action of the symmetric group 
$S_{[{n\over 2}]}$ on Tab$_n^{\bf D}$. This action is parallel 
to the action of $S_n$ on the ordinary tableaux studied in [2],
and has a natural interpretation in terms of self--evacuating tableaux
(Theorem~1.8).  One can conjecture that for even $n$ our action
coincides with the Carre--Leclerc's action.

The material of the paper is arranged as follows. In Section 1 we list main 
results and construct the bijection.  

Section 2 contains proofs of Theorems~1.6 and 1.8 on the action of  
the symmetric group $S_{[{n\over 2}]}$. 

For the reader's convenience, in Appendix we remind definitions of the
Bender--Knuth's involution and domino tableaux. 
\bigskip 

{\bf Acknowledgements}. The authors are grateful to George Lusztig for 
fruitful discussions. Special thanks are due to Mark Shimozono for an 
inestimable help  in the preparation of the manuscript. Both authors are 
grateful to CRM and the Department of Mathematics of University of Montreal, 
for hospitality and support during this work. 
\vfil\eject
 
\sssectitle {1. Main results}


Following [10], denote by $Tab_n$ the set of all (semi-standard Young) 
tableaux with entries not exceeding $n$.  For a partition $\lambda$ with at 
most $n$ rows, and for an integer vector $\beta=(\beta_1,\ldots,\beta_n)$ 
denote by  $Tab_\lambda(\beta)$ the set of $T\in Tab_n$ of  shape $\lambda$ 
and  weight $\beta$ (for the reader's convenience we collect all necessary 
definitions in the Appendix).  

Let $t_i:Tab_n\to Tab_n$, ($i=1,\ldots,n-1$) be the {\it Bender-Knuth} 
automorphisms which we define in the Appendix (see also [1],[2] and [6]). 

Define the automorphism ${\bf D}:Tab_n\to Tab_n$ by the formula:  
$${\bf D}={\bf D}_n:=t_{n-1}t_{n-3}\cdots \ , \eqno (1.1)$$
i.e., ${\bf D}=t_{n-1}t_{n-3}\cdots t_2$ for odd $n$, and ${\bf D}=
t_{n-1}t_{n-3}\cdots t_1$ for even $n$.

Denote by of  $Tab_n^{\bf D}$ the set of tableaux $T\in Tab_n$ 
such that ${\bf D}(T)=T$. For every shape $\lambda$ and weight $\beta'=
(\beta_1,\ldots,\beta_{[{n+1\over 2}]})$ denote $Tab_\lambda^{\bf D}(\beta')
=Tab_\lambda(\beta)\cap Tab_n^{\bf D}$, where 
$\beta=(\beta_{[{n+1\over 2}]},\cdots,\beta_k,\beta_k,\cdots,\beta_2,
\beta_2,\beta_1,\beta_1)$.

The {\it semi-standard domino tableaux} of shape $\lambda$ and weight 
$\beta'$ were studied by several autors (see, e.g., [10], p. 139). 
We also define them in the Appendix.  Our first "result" is another definition 
of domino tableaux.

\proclaim Proposition-Definition 1.1. For every $\lambda$ and $\beta'$ 
as above the domino tableaux of  shape $\lambda$ and weight $\beta'$ are in a 
natural one-to-one correspondence with the set $Tab_\lambda^{\bf D}(\beta')$. 

\medskip 

We prove this proposition in the Appendix. From now on we identify the 
domino tableaux with the elements of $Tab_n^{\bf D}$. 

Let us introduce the last bit of notation prior to the next result.

Let $p_i:=t_1t_2\cdots t_i$ for 
$i=1,\ldots,n-1$ (this $p_i$ is the inverse of the 
$i$-th {\it promotion operator} defined in [11]). 
Define the automorphism ${\bf P}:Tab_n\to Tab_n$ by the formula:  
$${\bf P}={\bf P}_n:=p_{n-1}p_{n-3}\cdots , \eqno (1.2)$$ 
so  
${\bf P}=p_{n-1}p_{n-3}\cdots p_2$ for odd $n$, and ${\bf P}=p_{n-1}p_{n-3}
\cdots p_1$ for even $n$. 

It is well-known (see, e.g., [2],[6]) that the Sch\"utzenberger evacuation  
involution factorizes as follows.
$${\bf S}={\bf S}_n=p_{n-1}p_{n-2}\cdots p_1 \eqno (1.3)$$
We regard (1.3) as a definition of ${\bf S}$. 

Denote by $Tab_n^{\bf S}\i Tab_n$ the set of the {\it  self-evacuating} 
tableaux, that is all the tableaux $T\in Tab_n$ fixed under ${\bf S}$. 
Similarly, denote $Tab_\lambda^{\bf S}(\beta)=Tab_\lambda(\beta)\cap 
Tab_n^{\bf S}$. It is well-known that this set is empty unless 
$\beta_i=\beta_{n+1-i}$ for  $i=1,2,\ldots,[{n+1\over 2}]$. 
 
\proclaim Theorem 1.2. The above automorphism ${\bf P}$ of $Tab_n^{\bf D}$ 
induces the bijection $$Tab_n^{\bf D}\cong Tab_n^{\bf S} \ .  \eqno (1.4)$$
More precisely, ${\bf P}$ induces the bijection 
$Tab_\lambda^{\bf D}(\beta')\cong Tab_\lambda^{\bf S}(\beta) $ 
for any shape $\lambda$ and the weight 
$\beta'=(\beta_1,\ldots,\beta_{[{n+1\over 2}]})$, where 
$\beta=(\beta_1,\ldots,\beta_{[{n+1\over 2}]},\beta_{[{n\over 2}]},
\ldots,\beta_1)$.

\bigskip 

{\bf Remark}. Theorem 1.2 is the "realization" of the theorem in Section 0.
 
The proof of this theorem is so elementary that we present it here. 

\medskip 

\noindent {\sl Proof of Theorem 1.2}. 

The statement (1.4) immediately follows from the surprisingly elementary 
lemma below. 

\proclaim Lemma 1.3. ${\bf S}={\bf P}{\bf D}{\bf P}^{-1}$. 

Indeed, if $T\in Tab_n^{\bf D}$ then $T={\bf D}(T)={\bf P}^{-1}{\bf S}
{\bf P}(T)$ by Lemma 1.3, so ${\bf S}({\bf P}(T))={\bf P}(T)$, that is,  
${\bf P}(T)\in Tab_n^{\bf S}$. This proves the inclusion  
${\bf P}(Tab_n^{\bf D})\i Tab_n^{\bf S}$. 
The opposite inclusion  also follows. 

{\bf Remark.} By definition (see the Appendix), the involutions 
$t_1,\ldots,t_{n-1}$ involved in Lemma 1.3, satisfy the 
following obvious relations  (see e.g. [2]): 
$$t_i^2={\rm id},~~t_it_j=t_jt_i \ , ~~~~~1\le i,j \le  n-1,~~|i-j|>1  
\ . \eqno (1.5) $$ 
In fact, the factorization in Lemma 1.3 is based only on this relations 
(as we can see from the proof below), and, hence, makes sence in any group 
with the relations (1.5).

It also follows from the definition of $t_i$ that  
$$t_i:Tab_\lambda(\beta)\cong Tab_\lambda((i,i+1)(\beta) \ .$$
This proves the second assertion of Theorem 1.2  provided that Lemma 1.3  
is proved.

\smallskip 

\noindent {\sl Proof of Lemma 1.3}.  
We will proceed by induction on $n$. If $n=2$ then ${\bf S}=t_1$ 
while ${\bf P}=t_1$, and ${\bf D}=t_1$. If $n=3$ then ${\bf S}_3=t_1t_2t_1$ 
while ${\bf D}=t_2$, and ${\bf P}=t_1$.  
Then, for any $n>3$ the formula (1.3) can be rewritten as 
${\bf S}_n=p_{n-1}{\bf S}_{n-1}$. It is easy to derive from (1.5) that 
${\bf S}_n={\bf S}_{n-1}p_{n-1}^{-1}$. Applying both of the above relations, 
we obtain 
$${\bf S}_n=p_{n-2}{\bf S}_{n-2}p_{n-1}^{-1} \ . \eqno (1.6)$$

Finally, by (1.5) and the inductive assumption for $n-2$, 
$${\bf S}_n=p_{n-2}{\bf S}_{n-2}p_{n-1}^{-1}=p_{n-2}{\bf P}_{n-2}
{\bf D}_{n-2}{\bf P}_{n-2}^{-1}p_{n-1}^{-1}
=p_{n-2}{\bf P}_{n-2}t_{n-1}{\bf D}_n{\bf P}_{n-2}^{-1}p_{n-1}^{-1}  
$$
since ${\bf D}_{n-2}=t_{n-1}{\bf D}_n$. Finally, ${\bf P}_{n-2}^{-1}
p_{n-1}^{-1}={\bf P}_n^{-1}$, and $p_{n-2}{\bf P}_{n-2}t_{n-1}=
p_{n-2}t_{n-1}{\bf P}_{n-2}=p_{n-1}{\bf P}_{n-2}={\bf P}_n$. 

Lemma 1.3 is proved. 

Theorem 1.2 is proved.

\bigskip 

Recall (see the Appendix) that the weight of a domino tableau 
$T\in Tab_n^{\bf D}$ 
is the vector  $\beta'=(\beta_1,\ldots, \beta_{[{n+1\over 2}]})$, 
where $\beta_1$ is 
the number of occurences of $n$, 
$\beta_2$ is that of $n-2$, and so on. Denote by $Tab_n^{\bf D}(\beta')$ 
the set 
of all $T\in Tab_n^{\bf D}$ of weight 
$\beta'$. Our next main  result is the following

\proclaim Theorem 1.4. There exists a natural faithful action of the 
symmetric group 
$S_{[{n\over 2}]}$ on $Tab_n^{\bf D}$ preserving the shape and acting on  
weight by permutation. 

In order to define the action precisely, let us recall one of the 
results of [2]. 

\proclaim Theorem 1.5 ([2]). There is an action of $S_n$ on $Tab_n$ 
given by $(i,i+1)\mapsto s_i$, where  
$$s_i={\bf S}_it_1{\bf S}_i \ , ~~~~~~i=1,\ldots,n \ , \eqno (1.7)$$
that is, the $s_i$ are involutions, and $(s_js_{j+1})^3={\rm id}, 
s_is_j=s_js_i$ for all $i,j$ such that $|j-i|>1$. 
(Here ${\bf S}_i$ given by (1.3) acts on $Tab_n$). 

\bigskip 

{\bf Definition}. For $n\ge 4$ define the automorphisms 
$\sigma_i:Tab_n\to Tab_n$, $i=2,\ldots,n-2$  by the formula
$$\sigma_i:=t_is_{i-1}s_{i+1}t_i \ .\eqno (1.8)$$

{\bf Remark}. The automorphisms $s_i$ were first defined by Lascoux and 
Sch\"utzenberger in [8] in the context of plactique monoid theory. The 
definition (1.7) of $s_i$ first appeared in [2] in a more general situation. 
As a matter of fact, for the restriction of $s_i$ to the set of tableaux, 
(1.7) is implied in [8].
 
Theorem 1.4 follows from a much stronger result, namely the following theorem:

\proclaim Theorem 1.6. \item {(a)}  
${\bf D}\sigma_i=\sigma_i{\bf D}$ if and only if 
$i\equiv n~({\rm mod}~2)$. 
\item {(b)}
The automorphisms $\sigma_2,\sigma_3,\cdots \sigma_{n-2}$ are involutions 
satisfying the Coxeter relations $(\sigma_i\sigma_j)^{n_{ij}}={\rm id}$
for $i,j=2,\ldots,n-2\ , $ where 
$$n_{ij}=\cases { 3~~~~{\rm if}~ |j-i|=1~ or ~2 \cr 
                  6~~~~{\rm if} ~|j-i|=3 \cr 
		2~~~~{\rm if}~ |j-i|>3
		} \ . \eqno (1.9)$$ 

We will prove Theorem 1.6 in Section 2. 

{\bf Remark}. One can conjecture that the relations (1.8) (together with 
$\sigma_i^2={\rm id}$) give a presentation of the group $\Sigma_n$ 
generated by the $\sigma_i$.

Indeed, the action of $S_{[{n\over 2}]}$ on $Tab_n$ from Theorem 1.4 can be 
defined by    
$$(i,i+1)\mapsto \sigma_{n-2i} \ , ~~~~~~i=2,\ldots,n-2 \ .$$  
According to Theorem 1.6(a), this action preserves $Tab_n^{\bf D}$. By 
definition (1.7), this action preserves shape of tableaux and permutes 
the weights as  follows: 
$$\sigma_{n-2i}:Tab_\lambda^{\bf D}(\beta')\cong Tab_\lambda^{\bf D}((i,i+1)
(\beta')) \ .\eqno (1.10)$$

\bigskip 

Let us try to get the action of $S_{[{n\over 2}]}$ on $Tab_n^{\bf S}$ using 
the following nice property of the involutions $s_1,\ldots,s_{n-1}$ (which 
can be found, e.g., in [2], Proposition 1.4):
$${\bf S}s_i{\bf S}=s_{n-i}\ , ~~~~~~i=1,\ldots,n-1 \ . \eqno (1.11)$$
Let $\tau_1,\ldots,\tau_{[{n\over 2}]-1}$ be defined by $\tau_k:=s_ks_{n-k}$.
The following proposition is obvious. 

\proclaim Proposition 1.7. \item {(a)} ${\bf S}\tau_k=\tau_k{\bf S}$ for all 
$k$. 
\item {(b)} The group generated by the $\tau_k$ is isomorphic to 
$S_{[{n\over 2}]}$ via 
$$(k,k+1)\mapsto \tau_k \ ,~~~~~~~ k=1,\ldots, {[{n\over 2}]}-1 \ . 
\eqno (1.12)$$

\medskip 

Restricting the action (1.12) to $Tab_n^{\bf S}$ we obtain the desirable 
action of $S_{[{n\over 2}]}$ on $Tab_n^{\bf S}$.  

Our last main result compares of the latter action with that on 
$Tab_n^{\bf D}$. The following theorem confirms  naturality of 
our choice of $\sigma_j$ in (1.8).  

\proclaim Theorem 1.8. For $k=1,\ldots,[{n\over 2}]-1$, we have 
${\bf P}^{-1}\tau_k{\bf P}=\sigma_{n-2k}$.

\sssectitle {2. Proof of Theorems 1.6, 1.8 and related results}


We start with a collection of formulas relating $s_i$ and $t_j$. First of all, 
$$s_it_j=t_js_i \eqno (2.1)$$
for all $i,j$ with $|j-i|>1$ (This fact follows from [8] 
and is implied in [2], Theorem~1.1). Then, by definition (1.8) of $\sigma_i$, 
$$\sigma_it_j=t_j\sigma_i \eqno (2.2)$$ 
whenever $|j-i|>2$. 
Furthermore,  (1.3) and (1.7) imply that 
$$s_j=p_js_{j-1}p_j^{-1} \ , ~~~~j=1,\ldots,n-1 \ . \eqno (2.3)$$ 
  More generally, 
$$s_i=p_js_{i-1}p_j^{-1} \ , ~~~~1\le i\le j <n \ . \eqno (2.4)$$  

Finally, (2.3) and (2.1) imply that  
$$t_{j-1}s_jt_{j-1}=t_js_{j-1}t_j \ . \eqno (2.5)$$

{\bf Proof of Theorem 1.6(a)}. 
Let $i$ be congruent $n$ modulo $2$. Then by (2.2) and the definition (1.1) 
of ${\bf D}$,  
$${\bf D}\sigma_i{\bf D}=t_{i+1}t_{i-1}\sigma_it_{i+1}t_{i-1} \ . \eqno (2.6)$$
We now prove  that (2.6) is equal to $\sigma_i$. Indeed, 
$${\bf D}\sigma_i{\bf D}=t_{i+1}t_{i-1}\sigma_i t_{i+1}t_{i-1}=
t_{i+1}t_{i-1}t_is_{i-1}s_{i+1}t_it_{i+1}t_{i-1} \ .$$
Now we will apply the relations (2.5) in the form  $t_{i-1}t_is_{i-1}=
s_it_{i-1}t_i$, and  $s_{i+1}t_it_{i+1}=t_it_{i+1}s_i$. 
We obtain 
$${\bf D}\sigma_i{\bf D}=t_{i+1}(s_it_{i-1}t_i)(t_it_{i+1}s_i)t_{i-1}
=t_{i+1}s_it_{i-1}t_{i+1}s_it_{i-1}=t_{i+1}s_it_{i+1}t_{i-1}s_it_{i-1} \ .$$
Applying (2.5) with $j=i-1,i$ once again, we finally obtain  
$${\bf D}\sigma_i{\bf D}=(t_is_{i+1}t_i)(t_is_{i-1}t_i)=t_is_{i+1}s_{i-1}t_i
=t_is_{i-1}s_{i+1}t_i=\sigma_i$$
by Theorem 1.5, which proves that ${\bf D}\sigma_i=\sigma_i{\bf D}$. 
The implication "if" of {\bf Theorem  1.6(a)} is proved.  
Let us prove the "only if" implication. Let $\rho$ be the canonical
homomorphism 
from the group generated by $t_1,\ldots,t_{n-1}$ to $S_n$ defined by 
$t_j\mapsto (j,j+1)$.  Then $\rho(\sigma_i)=(i,i+1)(i-1,i)(i+1,i+2)(i,i+1)$ 
for $i=2,\ldots,n-2$, and $\rho({\bf D})=(n-1,n)(n-3,n-2) \cdots~$.  
If $i\not \equiv n~(mod ~2)$ then $\rho(\sigma_i)$ and $\rho({\bf D})$ do 
not commute. Hence $\sigma_i{\bf D}\ne {\bf D}\sigma_i$ for such $i$. Thus 
{\bf Theorem  1.6(a)} is proved.

\smallskip 

{\bf Proof of Theorem  1.6(b)}. We are going to consider several cases when  
$j-i=1,2,3$ and $>3$.
In each case we compute the product $\sigma_i\sigma_j$ up to conjugation, 
and express it as a product of several $s_k$. These products are always of 
finite order, because, according to Theorem 1.5, these $s_k$'s generate the 
group isomorphic to $S_n$. 

The easiest case is $j-i>3$. Then all the ingredients of $\sigma_j$ in (1.7) 
(namely $t_j,s_{j-1},s_{j+1}$) commute with those of $\sigma_i$ 
(namely $t_i,s_{i-1},s_{i+1}$). Thus $\sigma_j$ commutes with $\sigma_i$, 
that is, $(\sigma_i\sigma_j)^2={\rm id}$, and we are done in this case. 

In what follows we will write $a\equiv b$ if $a$ is conjugate to $b$.  

{\bf Case} $j-i=3$, i.e., $j=i+3$. We have
$$\sigma_i\sigma_{i+3}
=t_is_{i-1}s_{i+1}t_it_{i+3}s_{i+2}s_{i+4}t_{i+3}\equiv s_{i-1}s_{i+1} 
s_{i+2}s_{i+4}$$
by the above commutation between remote $s_k$ and $t_l$.

By Theorem 1.5, 
$(\sigma_i\sigma_{i+3})^6\equiv (s_{i-1})^6(s_{i+1} s_{i+2})^6(s_{i+4})^6
={\rm id}$.
So $(\sigma_i\sigma_{i+3})^6={\rm id}$ and we are done with this case. 

{\bf Case} $j-i=2$, i.e., $j=i+2$.  We have
$$\sigma_i\sigma_{i+2}
=t_is_{i-1}s_{i+1}t_it_{i+2}s_{i+1}s_{i+3}t_{i+2}\equiv t_{i+2}t_is_{i-1}
s_{i+1}t_it_{i+2}s_{i+1}s_{i+3}$$
$$=t_is_{i-1}t_{i+2}s_{i+1}t_{i+2}t_is_{i+1}s_{i+3}
=t_is_{i-1}t_{i+1}s_{i+2}t_{i+1}t_is_{i+1}s_{i+3}$$
by (2.5) with $j=i+2$. Conjugating the latter 
expression with $t_i$, we obtain, 
$$\sigma_i\sigma_{i+2} \equiv s_{i-1}t_{i+1}s_{i+2}t_{i+1}t_is_{i+1}s_{i+3} 
t_i=s_{i-1}t_{i+1}s_{i+2}t_{i+1}t_is_{i+1}t_is_{i+3}$$
$$=s_{i-1}t_{i+1}s_{i+2}t_{i+1}t_{i+1}s_it_{i+1}s_{i+3}
=s_{i-1}t_{i+1}s_{i+2}s_it_{i+1}s_{i+3}$$
$$\equiv t_{i+1}s_{i-1}t_{i+1}s_{i+2}s_it_{i+1}s_{i+3}t_{i+1}= s_{i-1}
s_{i+2}s_is_{i+3}=s_{i-1}s_is_{i+2}s_{i+3}  \ .$$
Thus, $(\sigma_i\sigma_{i+2})^3\equiv (s_{i-1}s_is_{i+2}s_{i+3})^3=
(s_{i-1}s_i)^3(s_{i+2}s_{i+3})^3={\rm id}$,
and we are done with this case too. 

{\bf Case} $j-i=1$, i.e., $j=i+1$. We have 
$$\sigma_i\sigma_{i+1}
=t_is_{i-1}s_{i+1}t_it_{i+1}s_is_{i+2}t_{i+1}
=t_is_{i-1}s_{i+1}t_it_{i+1}s_it_{i+1}t_{i+1}s_{i+2}t_{i+1}$$
$$=t_is_{i-1}s_{i+1}t_it_is_{i+1}t_it_{i+2}s_{i+1}t_{i+2}
=t_is_{i-1}t_it_{i+2}s_{i+1}t_{i+2}$$
$$=t_{i-1}s_it_{i-1}t_{i+2}s_{i+1}t_{i+2}\equiv 
s_it_{i+2}s_{i+1}t_{i+2}\equiv s_i s_{i+1} \ . $$ 
Thus, $(\sigma_i\sigma_{i+1})^3\equiv (s_is_{i+1})^3={\rm id}$, 
and we are done with this final case. 

This finishes the proof of {\bf Theorem  1.6}. 

\bigskip 

{\bf Proof of Theorem 1.8}. First of all we prove Theorem 1.8 with $k=1$, 
that is, 
$${\bf P}_n^{-1}s_1s_{n-1}{\bf P}_n=\sigma_{n-2} \eqno (2.7)$$
Since ${\bf P}_n=p_n{\bf P}_{n-2}$, the left hand side of (2.7) equals 
$${\bf P}_{n-2}^{-1}p_{n-1}^{-1}s_1p_{n-1}p_{n-1}^{-1}s_{n-1}p_{n-1}
{\bf P}_{n-2}={\bf P}_{n-2}^{-1}(p_{n-1}^{-1}s_1p_{n-1})s_{n-2}
{\bf P}_{n-2} \ . \eqno (2.8)$$
Now we prove by induction on $n$ that 
$$p_{n-1}^{-1}s_1p_{n-1}=p_{n-2}s_{n-1}p_{n-2}^{-1} \eqno (2.9)$$
Indeed, if $n=2$, (2.9) is obvious (we agree that $p_0=t_0={\rm id}$). 
Then using (2.9) as inductive assumption, we obtain 
$$p_n^{-1}s_1p_n=t_n(p_{n-1}^{-1}s_1p_{n-1})t_n
=t_n(p_{n-2}s_{n-1}p_{n-2}^{-1})t_n$$
$$=p_{n-2}t_ns_{n-1}t_np_{n-2}^{-1}
=p_{n-2}t_{n-1}s_nt_{n-1}p_{n-2}^{-1}=p_{n-1}s_np_{n-1}^{-1}$$
by (2.5), and we are done with (2.9). 

Finally, in order to finish the proof of (2.7), let us substitute (2.9) 
into (2.8).  Thus, taking into account that ${\bf P}_{n-2}=p_{n-3}
{\bf P}_{n-4}=p_{n-2}t_{n-2}{\bf P}_{n-2}$, (2.8) is equal to 
$${\bf P}_{n-2}^{-1}(p_{n-2}s_{n-1}p_{n-2}^{-1})s_{n-2}{\bf P}_{n-2}=
{\bf P}_{n-4}^{-1}t_{n-2}p_{n-2}^{-1}(p_{n-2}s_{n-1}
p_{n-2}^{-1})s_{n-2}p_{n-2}t_{n-2}{\bf P}_{n-4} $$
$$={\bf P}_{n-4}^{-1}t_{n-2}s_{n-1}
(p_{n-2}^{-1}s_{n-2}p_{n-2})t_{n-2}{\bf P}_{n-4}={\bf P}_{n-4}^{-1}t_{n-2}
s_{n-1}s_{n-3}t_{n-2}{\bf P}_{n-4}$$
$$=t_{n-2}s_{n-1}s_{n-3}t_{n-2}=\sigma_{n-2} \ , $$
which proves {\bf Theorem 1.8} for $k=1$. 

Now using induction on $n$ we will prove the general case, namely, the formula 
$${\bf P}_n^{-1}s_ks_{n-k}{\bf P}_n=\sigma_{n-2k} \ , 1\le k \le 
[{n\over 2}]-1 . \eqno (2.10)$$
Indeed, 
$${\bf P}_n^{-1}s_ks_{n-k}{\bf P}_n={\bf P}_{n-2}^{-1}(p_{n-1}^{-1}s_k
s_{n-k}p_{n-1}){\bf P}_{n-2} . \eqno (2.11)$$
But $p_{n-1}^{-1}s_ks_{n-k}p_{n-1}=(p_{n-1}^{-1}s_kp_{n-1})(p_{n-1}^{-1}
s_{n-k}p_{n-1})=s_{k-1}s_{n-k-1}$ by (2.4) applied with $j:=n-1, i=k$, and 
$i=n-k$. 
Thus (2.11) is equal to 
$${\bf P}_{n-2}^{-1}s_{k-1}s_{n-k-1}{\bf P}_{n-2}=\sigma_{n-2-2(k-1)}=
\sigma_{n-2k}$$
by the inductive assumption (2.10) with $n-2$. 
{\bf Theorem 1.8} is proved. 


\sssectitle {3. Appendix. Definitions of domino tableaux} 


By definition from [10], a {\it (semi-standard  Young) tableau} $T$ {\it with
entries not exceeding} $n$ is an ascending sequence of Young diagrams
$T=(\{\emptyset\}=\Lambda_0\i \Lambda_1\i \ldots \i \Lambda_n)$ such 
that $\Delta_i=\Lambda_i\setminus \Lambda_{i-1}$ is a {\it
horizontal strip} for $i=1,\ldots,n$ (we allow
$\Lambda_{i+1}=\Lambda_i$). Denote by $Tab_n$ the set of all
(semi-standard Young) tableaux with entries not exceeding $n$. 
Define the {\it shape} of $T$ by $\lambda=\Lambda_n$, and the {\it weight} 
$\beta=(\beta_1,\beta_2,\ldots,\beta_n)$ of $T$ by 
$\beta_i=|\Delta_i|$ for $i=1,\ldots,n$. Denote by $Tab_n(\beta)$ 
the set of all
$T\in Tab_n$ having weight $\beta$. 
 
One visualizes the tableau $T$  as filling of each box of
$\Delta_i$ with the number $i$ for $i=1,\ldots,n$.

Let $T=(\Lambda_0\i  \Lambda_1 \i\Lambda_2  )$ be a {\it skew} 
$2$-{\it tableau}, that is, 
$\Lambda_1\setminus \Lambda_0$ and $\Lambda_2\setminus \Lambda_1$ both are 
horizontal strips. 
For such a skew $2$-tableau $T$ define the skew $2$-tableau
$t(T)=T=(\Lambda_0\i \Lambda'_1\i \Lambda_2)$ as follows. 
First of all we visualize $T$ as the shape $D=\Lambda_2\setminus 
\Lambda_0$ filled
with the letters ${\bf a}$ 
and ${\bf b}$ where the letters ${\bf a}$ fill  the horizontal strip
$\Lambda_1\setminus \Lambda_0$, and the letters ${\bf b}$ fill the
remaining horizontal strip $\Lambda_2\setminus \Lambda_1$. 
Then $t(T)$ will be a new filling of $D$ with the same letters ${\bf a}$ 
and ${\bf b}$ which we construct below.
 
Let $D_k$ be the longest connected sub-row of the $k$-th  row of $D$,
whose  boxes have no horizontal edges in common with other boxes of
$D$. $D_k$ contains say, $l$ boxes on the left filled with ${\bf a}$, and
$r$ boxes (on the right) 
filled with  ${\bf b}$. We will replace this filling of $D_k$ with the
$r$ letters ${\bf a}$ on the left, and $l$ letters ${\bf b}$ on the
right. 
We do this procedure for every $k$, and leave unchanged the filling of the 
complement of the union of all the $D_k$ in $D$. 
It is easy to see that the new filling
of $D$ is a skew $2$-tableau of the form $t(T)=T=(\Lambda_0\i \Lambda'_1\i 
\Lambda_2)$. 
Clearly the correspondence $T\mapsto t(T)$ is an involutive automorphism
of the set of skew $2$-tableaux. 

The following lemma is obvious. 

\proclaim Lemma-Definition A1. Let $D$ be a skew Young diagram.  Then the
following are equivalent:
\item {(i)} There is a (unique)  skew $2$-tableau 
$T=(\Lambda_0\i \Lambda_1\i \Lambda_2)$ with $D=\Lambda_2\setminus 
\Lambda_0$ such that 
$$t(T)=T$$
\item {(ii)} There is a (unique) covering of $D$ by (non-overlapping) dominoes 
such that every domino has no common horizontal edges with other
boxes of $D$.

Define the {Bender-Knuth involution} ([1], see also [2]) 
$t_i:Tab_n\to Tab_n$ as follows.  
 $T=(\cdots \i \Lambda_{i-1}\i \Lambda_i \i
\Lambda_{i+1} \i \cdots )\in Tab_n$. Define  
$T_i(T)\in Tab_n$ by 
$$t_i(T)=(\cdots \i \Lambda_{i-1}\i \Lambda'_i \i
\Lambda_{i+1}\i \cdots )$$
where $\Lambda_i'$ is defined by $(\Lambda_{i-1}\i \Lambda'_i\i
\Lambda_{i+1})=t(\Lambda_{i-1}\i \Lambda_i \i
\Lambda_{i+1})$, and other $\Lambda_j$ remain unchanged.

Denote the resulting tableau by $t_i(T)$. 
It is easy to see that the correspondence $T\mapsto t_i(T)$ is a
well-defined involutive automorphism $t_i:Tab_n\to Tab_n$, and
$t_it_j=t_jt_i$ whenever $|j-i|>1$.

 A {\it domino tableau}
${\bf T}$ {with entries not exceeding} ${n+1\over 2}$ is an ascending   
sequence of the Young  diagrams ${\bf T}=(\Lambda_{\varepsilon}\i 
\Lambda_{\varepsilon+2} \i \cdots \i \Lambda_{n-2}\i  \Lambda_n)$ 
with $\varepsilon=n-2[{n\over 2}]$, such that

\item {(i)}  $\Lambda_{\varepsilon}=\emptyset$ if $n$ is even, 
(i.e., $\varepsilon=0$), and $\Lambda_{\varepsilon}$ is a connected horizontal
strip  if $n$ is odd, (i.e., $\varepsilon=1$);

\item {(ii)} Each difference $D_k=\Lambda_{n-2k}\setminus \Lambda_{n-2k+2}$ 
is subject to
condition (ii) of Lemma A1. 

\bigskip 

The {\it shape} of ${\bf T}$ is defined to be $\Lambda_n$. We also define 
the {\it weight} $\beta'=(\beta_1,\ldots,\beta_{[{n+1\over 2}]})$ of 
${\bf T}$ by  $\beta_k={|D_k|\over 2}$ for $k=1,\ldots,[{n+1\over 2}]$.    

Thus the above definitions together with Lemma-Definition A1 make
Proposition-Definition 1.1 obvious.

\vfil\eject
 
\sssectitle { References} 

\bigskip 

\bitem {1}
E.~Bender, D.~Knuth, Enumeration of plane partitions. {\sl J. 
Combinatorial Theory} Ser. A, {\bf 13} (1972), 40--54.

\bitem {2} 
A.~Berenstein, A.N.~Kirillov, Groups generated by involutions,
Gelfand-Tsetlin patterns and combinatorics of Young tableaux,  
{\sl Algebra i Analiz}, {\bf 7} (1995), no. 1, 92--152; translation in St. 
Petersburg Math. J. {\bf 7} (1996),
no. 1, 77--127.

\bitem {3} 
A.~Berenstein, A.~Zelevinsky,  Canonical bases for the quantum group of 
type $A\sb r$ and piecewise-linear
combinatorics, {\sl Duke Math. J.}, {\bf 82} (1996), no. 3, 473--502.

\bitem {4}
C.~Carre, B.~Leclerc, Splitting the square of Schur function into
symmetric and antisymmetric parts, {\sl J. Algebraic Comb.}, {\bf 4}
(1995), n.3, 201--231.

\bitem  {5}
S.~Fomin, Generalized Robinson-Schensted-Knuth correspondence, 
{\sl Journal of Soviet Math.}, {\bf 41} (1988), 979--991.

\bitem  {6} 
E.~Gansner, On the equality of two plane partition
correspondences, {\sl Discrete Math.}, {\bf 30} (1980), 121--132.

\bitem  {7}
A.N.~Kirillov, A.~Lascoux, B.~Leclerc, J.-Y.~Thibon, Series 
generatrices pour les tab\-leaux de domino, {\sl C.R.Acad. Sci. Paris},
{\bf 318} (1994), 395--400.

\bitem {8} 
A.~Lascoux, M.-P.~Sch\"utzenberger,  Le monoïde plaxique. 
{\sl Noncommutative structures in algebra and geometric combinatorics} 
(Naples, 1978), pp. 129--156, Quad. "Ricerca Sci.", 109, CNR, Rome, 1981.

\bitem {9}
M.A. van Leeuwen, The Robinson-Schensted and Sch\"utzenberger algorithms,
an elementary approach, {\sl Electronic Journal of Combinatorics}, {\bf 3},
n.2 (1996).
 
\bitem {10} 
I.~Macdonald, {\it Symmetric functions and Hall polynomials}. Second edition. 
With contributions by A.
Zelevinsky. Oxford Mathematical Monographs. Oxford Science Publications. 
The Clarendon Press, Oxford University
Press, New York, 1995. 
 
\bitem  {11}
M.-P.~Sch\"utzenberger, Promotion des morphismes
d'ensembles ordonn\'es, {\sl Discrete \hfill \break Math.}, 
{\bf 2} (1972), 73--94.
 
\bitem  {12} 
J.~Stembridge, Canonical bases and self-evacuating tableaux, 
{\sl Duke Math. J.}, {\bf 82} (1996), no.3, 585--606.

\end